\documentclass[11pt]{amsart}
\usepackage{amssymb,latexsym}

\title{Kodaira dimension in low dimensional  topology}

\newcommand{\C}{\mathbb{C}}

\newcommand{\Q}{\mathbb{Q}}
\newcommand{\R}{\mathbb{R}}
\newcommand{\Z}{\mathbb{Z}}

\def\be{\begin{equation}}
\def\ee{\end{equation}}

\newtheorem{theorem}{Theorem}[section]

\newtheorem{lemma}[theorem]{Lemma}
\newtheorem{remark}[theorem]{Remark}
\newtheorem{question}[theorem]{Question}
\newtheorem{definition}[theorem]{Definition}

\newtheorem{conj}[theorem]{Conjecture}

\begin{document}

\author{ Tian-Jun Li}
\address{School  of Mathematics\\  University of Minnesota\\ Minneapolis, MN 55455}
\email{tjli@math.umn.edu}

\date{\today}
\maketitle

\begin{abstract}
This is a survey on the various notions of Kodaira dimension in low dimensional topology. 
The focus is on  progress  after the 2006 survey \cite{survey:Kod-dim4}.  
\end{abstract}

\tableofcontents

\section {Introduction}
Roughly speaking, a Kodaira dimension type invariant on a class of $n-$dimensional manifolds
is a numerical invariant taking  values in the finite set

$$\{-\infty,   0, 1,  \cdots,   \lfloor\frac{n}{2}\rfloor\},$$
where $\lfloor x\rfloor$ is
the largest integer bounded by $x$.

The first invariant of this type is due to Kodaira (\cite{K}) for smooth algebraic varieties  (naturally extended to complex manifolds): Suppose $(M, J)$ is a complex manifold of real dimension $2m$.
The holomorphic Kodaira dimension $\kappa^{h}(M,J)$ is defined as
follows:

\[
\kappa^{h}(M,J)=\left\{\begin{array}{ll}
-\infty &\hbox{ if $P_l(M,J)=0$ for all $l\ge 1$},\\
0& \hbox{ if $P_l(M,J)\in \{0,1\}$, but $\not\equiv 0$ for all $l\ge 1$},\\
k& \hbox{ if $P_l(M,J)\sim cl^k$; $c>0$}.\\
\end{array}\right.
\]
Here  $P_l(M,J)$ is the $l-$th plurigenus of the  complex manifold
$(M, J)$ defined by $P_l(M,J)=h^0({\mathcal K}_J^{\otimes l})$, with
${\mathcal K}_J$ the canonical bundle of $(M,J)$.

This classical Kodaira dimension $\kappa^h$, along with its various extensions, has played an essential role in algebraic geometry, K\"ahler geometry and complex geometry (cf. \cite{I}, \cite{De}, \cite{Ue}).

In the past 20 years, it has been gradually realized that such   notions also exist 
 in low dimensional topology. 
The  paper \cite{survey:Kod-dim4} is a survey on  the Kodaira dimension $\kappa^s$ for
symplectic 4-manifolds up to 2005.
The current article updates the progress on $\kappa^s$ and its extensions, as well as its cousin $\kappa^t$ for 3-manifolds. 

The construction of the symplectic Kodaira dimension $\kappa^s$  is impossible without Taubes' fundamental works. We would like to 
take this opportunity to express our deep  gratitude towards him. 
We also benefit from discussions with I. Baykur,  J. Dorfmeister, J. Fine, S. Friedl, Y. Koichi, G. LaNave, C. Mak, D. Salamon, W. Wu, W. Zhang. We are also grateful to the referee for useful suggestions. 
This work is supported by NSF. 

\section{$\kappa^t$ and $\kappa^s$}

Let $M$ be a closed, smooth, oriented manifold.
  To begin with, we
make the following definition for  logical compatibility.

\begin{definition}
If $M=\emptyset$, then its Kodaira dimension   is defined to be
$-\infty$.
\end{definition}

\subsection{The topological Kod dim $\kappa^t$ for manifolds up to dimension $3$}
\subsubsection{$\kappa^t$ in dimensions $0, 1$ and $2$}

 The only closed
connected $0-$dimensional manifold is a point, and the only closed
connected $1-$dimensional manifold is a circle.

\begin{definition} \label{01} If $M$ has dimension $0$ or $1$, then
its  Kodaira dimension $\kappa^t(M)$ is defined to be $0$.
\end{definition}

The $2-$dimensional Kodaira dimension is defined by the positivity
of the Euler class. 
Suppose $M^2$ is a $2-$dimensional closed, connected, oriented manifold with Euler class
$e(M^2)$. Write $K=-e(M^2)$ and define
\[ \kappa^t(M^2)=\left\{\begin{array}{cc}
-\infty & \hbox{ if $K<0$},\\
0& \hbox{ if $K=0$},\\
1& \hbox{ if $K> 0$}.
\end{array}\right.
\]

It is easy to see that for any complex structure $J$ on $M^2$, $K$
is its canonical class, and  $\kappa^h(M^2, J)=\kappa^t(M^2)$.
 $\kappa^t(M^2)$ can be further interpreted from other viewpoints:
 symplectic structure
($K$ is also the symplectic canonical class),
 the Yamabe invariant,
geometric structures etc.

\subsubsection{$\kappa^t$  in dimension 3}
We move on to dimension $3$. In this dimension the definition of the
Kodaira dimension in \cite{Z} by Weiyi Zhang  is based on
geometric structures in the sense of Thurston. 

Divide the $8$
Thurston geometries into $3$ categories: 
$$\begin{array}{ll}
    -\infty: &\hbox{$S^3$ and $S^2 \times \mathbb{R}$};\cr
     0: &\hbox{$\mathbb{E}^3$, Nil and Sol};\cr
     1: &\hbox{$\mathbb{H}^2\times \mathbb{R}$, $\widetilde{SL_2(\mathbb{R})}$ and
    $\mathbb{H}^3$}.
\end{array}$$
Given a closed, connected $3-$manifold $M^3$,  we decompose it first by a prime
decomposition and then further consider  a  toroidal decomposition
for each prime summand, such that at the end each piece has a
geometric structure either in group $(1)$, $(2)$ or $(3)$ with
finite volume. 

\begin{definition} For a closed, connected $3-$dimensional manifold $M^3$,  its Kodaira dimension $\kappa^t(M)$ is defined  as follows:

\begin{enumerate}
    \item $\kappa^t(M^3)=-\infty$ if for any decomposition, each piece has  geometry type
    in category $-\infty$,
    \item $\kappa^t(M^3)=0$ if for any decomposition, we have at least one piece with geometry type
    in category $0$, but no piece has type in category $1$,
    \item $\kappa^t(M^3)=1$ if for any decomposition, we have at
    least one piece in category $1$.
\end{enumerate}

\end{definition}

The following are basic properties and facts established in \cite{Z}:
\begin{itemize}

\item $\kappa^t$ is  additive for any fiber bundle. 

\item  If there is a nonzero degree map from $M$ to $N$, then $\kappa^t(M)\geq \kappa^t(N)$.

\item If $\kappa^t(M)=-\infty$, then each prime summand of $M$ is  either
(a) spherical, i.e. it has a Riemannian metric of constant positive sectional curvature;
or (b) an $S^2$  bundle over $S^1$. 

\item If $\kappa^t(M)=0$, then each prime summand of $M$ is either in (a) or (b), or it is 
(c) a Seifert fibration with zero orbifold Euler characteristic;
or (d) a mapping torus of an Anosov map of the 2-torus or quotient of
these by groups of order at most 8.

\end{itemize}

Let $vb_1(M)$ be the supremum of $b_1(\tilde M)$ among all finite covers $\tilde M$. 
Due to  Agol's remarkable resolution of the virtual Betti number conjecture, $\kappa^t$ has the following interpretation
in terms of the virtual 1st Betti number $vb_1$,  at least for irreducible 3-manifolds (private communication with W. Zhang):
$\kappa^t = -\infty $ when $vb_1 =0$, 
$\kappa^t=0$  when $vb_1$ is finite and positive,
$\kappa^t=1$ when $vb_1$  is infinite.

{Notice that  we use $\kappa^t$ to denote the Kodaira dimension
for smooth manifolds in dimensions $0, 1, 2, 3$. 
Here $t$ stands for
{\it topological/smooth}, because in these dimensions homeomorphic
manifolds are actually diffeomorphic.

\medskip
For a non-orientable, connected manifold, we define its Kodaira dimension to 
be that of its (unique) orientable, connected covering. 
For a possibly disconnected manifold, we define its
Kodaira dimension to be the maximum of that of its components.

\bigskip
 In
summary,  the Kodaira dimension is defined for all  closed smooth
manifolds with dimension less than $4$, and is a topological/smooth invariant.

\subsection{The symplectic Kod dim $\kappa^s$ for $4-$manifolds}
\subsubsection{Definition based on Taubes SW}

Let $M$ be a closed, oriented smooth 4-manifold. 
Let ${\mathcal E}_M$ be the set of cohomology
classes whose Poincar\'e dual are represented by smoothly embedded
spheres of self-intersection $-1$. $M$ is said to be (smoothly)
minimal if ${\mathcal E}_M$ is the empty set.
Equivalently, $M$ is minimal if it is not the connected sum of
another manifold  with $\overline{\mathbb {CP}^2}$.

\medskip
Suppose $\omega$ is a symplectic form compatible with the orientation. 
$(M,\omega)$  is said to be (symplectically) minimal if ${\mathcal
E}_{\omega}$ is empty, where
$${\mathcal E}_{\omega}=\{E\in {\mathcal
E}_M|\hbox{ $E$ is represented by an embedded $\omega-$symplectic
sphere}\}.$$ 
We say that $(N, \tau)$
is a minimal model of $(M, \omega)$ if $(N, \tau)$ is  minimal and $(M, \omega)$ is 
a symplectic blow up of $(N, \tau)$. 
A basic fact proved using Taubes SW theory is:
  ${\mathcal E}_{\omega}$ is
empty if and only if ${\mathcal E}_M$ is empty. In other words, $(M,
\omega)$ is symplectically minimal if and only if $M$ is smoothly
minimal.

For a minimal symplectic $4-$manifold $(M^4,\omega)$ with symplectic
canonical class $K_{\omega}$,   the Kodaira dimension of
$(M^4,\omega)$
 is defined in the following way:

$$
\kappa^s(M^4,\omega)=\begin{cases} \begin{array}{lll}
-\infty & \hbox{ if $K_{\omega}\cdot [\omega]<0$ or} & K_{\omega}\cdot K_{\omega}<0,\\
0& \hbox{ if $K_{\omega}\cdot [\omega]=0$ and} & K_{\omega}\cdot K_{\omega}=0,\\
1& \hbox{ if $K_{\omega}\cdot [\omega]> 0$ and} & K_{\omega}\cdot K_{\omega}=0,\\
2& \hbox{ if $K_{\omega}\cdot [\omega]>0$ and} & K_{\omega}\cdot K_{\omega}>0.\\
\end{array}
\end{cases}
$$
Here $K_{\omega}$ is defined as the first Chern class of the
cotangent bundle for any almost complex structure compatible with
$\omega$.

$\kappa^s$ is well defined  since there does not
exist a minimal $(M, \omega)$ with $$K_{\omega}\cdot [\omega]=0, 
\quad\hbox{and}\quad  K_{\omega}\cdot K_{\omega}>0.$$
This again follows from Taubes SW theory. 
Moreover, 
 $\kappa^s$ is independent of $\omega$, so it is an oriented diffeomorphism invariant of $M$. 
And it follows from \cite{Liu} (cf. also \cite{OO}) that  $\kappa^s(M)=-\infty$ if and only if  $M$ is $\mathbb {CP}^2, S^2\times S^2$ or an $S^2-$bundle over a Riemann surface
of positive genus.

The Kodaira dimension of a non-minimal manifold is defined to be
that of any of its minimal models.
$\kappa^s(M, \omega)$ is well-defined for any $(M, \omega)$ since 
 minimal models always  exist. Moreover,  minimal models are almost unique up to diffeomorphisms:   If $(M, \omega)$ has non-diffeomorphic minimal models, then these minimal models have $\kappa^s=-\infty$. 
Diffeomorphic  minimal models have the same $\kappa^s$. 



Here are basic properties of $\kappa^s$:

\begin{itemize}

\item $\kappa^s$  is an oriented diffeomorphism invariant of $M$ (\cite{L06}). 

\item  $\kappa^s=\kappa^h$ whenever both
are defined (\cite{DZ}).

\end{itemize}

We remark that it was shown  by Friedman and Qin (\cite{FQ})  that
$\kappa^h(M^4, J)$ only depends on the
oriented diffeomorphism type of $M^4$. In light of these properties of $\kappa^s$ and $\kappa^h$, we ask 

\begin{question} To what extent  can  $\kappa^s$ and
$\kappa^h$ be extended to $\kappa^d$ for  smooth $4-$manifolds (here $d$
stands for `differentiable').
\end{question}

\subsubsection{Yamabe invariant}

Recall that the Yamabe invariant is defined in the following way: 
\begin{equation}\label{yamabeinv}Y(M) =
\sup_{[g]\in \mathcal C}\inf_{g\in [g]}\int_M s_g
dV_g,\end{equation} where $g$ is a Riemannian  metric on $M$, $s_g$
the scalar curvature, $[g]$ the conformal class of $g$, and
$\mathcal C$ the set of conformal classes on $M$.

A basic fact is  that $Y(M)>0$ if and only if $M$ admits a metric of
positive scalar curvature. 
Thus  $Y(M)$ is non-positive if $M$ does
not admit metrics of positive scalar curvature. 
Furthermore, in this
case, another basic fact is that $Y(M)$ is the supremum of the
scalar curvatures of all unit volume constant-scalar-curvature
metrics on $M$ (such metrics exist due to the resolution of the
Yamabe conjecture). 
It immediately follows that, in dimension two,
the sign of $Y(M^2)$ completely determines $\kappa^t(M^2)$.

In dimension three, $Y(M^3)$ is also closely related to the   geometric
structure of $M^3$, at least when $M^3$ is irreducible.
However,  the number $Y(M^3)$ does not
completely determine $\kappa^t(M^3)$ (See \cite{Z} and \cite{LZ}). 

 When $M^4$ admits a  K\"ahler
structure, 
LeBrun (\cite{LeB2})
 calculated $Y(M^4)$ and concluded that \eqref{yamabeinv} completely
determines $\kappa^h$. As $\kappa^s=\kappa^h$ for a K\"ahler
surface, if $M^4$ admits a K\"ahler structure,
 then
\begin{equation}\label{lebrun}  \kappa^s(M^4)=\left\{\begin{array}{ll}
-\infty & \hbox{ if $Y(M^4)>0$},\\
0& \hbox{ if $Y(M^4)=0$ and $0$ is attainable by a metric},\\
1& \hbox{ if $Y(M^4)=0$ and $0$ is not attainable},\\
2& \hbox{ if $Y(M^4)<0$}.
\end{array}\right.
\end{equation}
It is worth mentioning that there is connection between Ricci flow and the Yamabe invariant,
which is related to the beautiful criterion \eqref{lebrun} in the K\"ahler context. 
(See eg. \cite{IRS}).
However, in the general symplectic context, \eqref{lebrun} does not determine $\kappa^s(M^4)$: All $T^2-$bundles over $T^2$ have $\kappa^s=0$, 
while most of them do not have any zero scalar curvature
metrics. 
And, while $\kappa^s$ is invariant under finite coverings,   the sign of the Yamabe invariant   is not a covering invariant
of 4-manifolds (\cite{LeB1}). 
But the question of LeBrun  still makes
sense:
 if $M^4$ admits a symplectic structure and $Y(M^4)<0$, is $\kappa^s(M^4)$ equal to $2$?
(cf. \cite{LeB}, \cite{Su}, \cite{Torres}). 


\subsection{Symplectic manifolds of dimension 6 and higher}

 In higher dimension,  Kodaira dimension is only defined for
complex manifolds. Further, it is known that  $\kappa^h$ is not a
diffeomorphism invariant  \cite{R}. 
Thus we can only expect to have a notion of Kodaira dimension for
smooth manifolds with some additional  structures. For higher dimensional symplectic manifolds, there is a proposal  to extend $\kappa^s$ in \cite{LR}.
Another approach via Donaldson' peak sections is investigated in 
\cite{LaL}.

\bigskip

In the rest of the paper we will focus on symplectic 4-manifolds.

\section{Calculating $\kappa^s$}
\subsection{Additivity and subadditivity}

\subsubsection
{Subadditivity for Lefschetz fibrations}

A central problem  in birational  geometry is the following  Iitaka conjecture $C_{n, m}$ for holomorphic fibrations of algebraic varieties: Let $f : X \to Z$  be an algebraic fibre space where $X$ and $Z$ are smooth projective varieties of dimension $n$ and $m$, respectively, and let $F$ be a general fibre of $f$. Then,
$\kappa^h(X) \geq \kappa^h(F) + \kappa^h(Z)$.

It has been verified  when $n\leq 6$  (cf.  \cite{B} for the status of this conjecture). 
For symplectic 4-manifolds, Lefschetz fibrations are the analogues of holomorphic fibrations, and it is established in \cite{DZ} that 
\begin{equation} \label{subadditivity}   \kappa^s(X)\geq \kappa^t(Base) +\kappa^t(Fiber).
\end{equation}

This is certainly true when the base is $S^2$. 
When the base genus is at least $1$, 
 given a relative minimal $(g, h, n)$ Lefschetz fibration with $h \geq 1$,  the Kodaira dimension $\kappa^l(g, h, n)$ is introduced in \cite{DZ}:
$$\kappa^l(g, h, n)=\begin{cases}  -\infty & \mbox{ if } g=0,\\
0 &\mbox{ if }  (g,h,n) = (1,1,0),\\
1 &\mbox{ if } (g,h) = (1, \geq  2)  \mbox{ or } (g,h,n) = (1,1,> 0) \mbox { or }  (\geq 2,1,0),\\
2 &\mbox{ if }   (g,h) \geq  (2,2)  \mbox{ or } (g,h,n) = (\geq 2,1, \geq  1).
\end{cases}
$$
The Kodaira dimension of a non-minimal Lefschetz fibration with $h\geq 1$
is defined to be that of its minimal models.
Further it is verified in \cite{DZ} and \cite{D-tori} that $\kappa^l(g, h, n)=\kappa^s$.
Clearly the subadditivity \eqref{subadditivity}   follows.

\subsubsection{Additivity for fibred manifolds}

For a  holomorphic
fiber bundle in the projective category, a classical theorem \cite{Ib} says that  the holomorphic Kodaira dimension $\kappa^h$
is additive.
The additivity for $\kappa^t$ has been established in \cite{Z}, and there is strong evidence that
it is also valid  for $\kappa^s$. 

\medskip
\noindent {\bf $0-$dimensional fibers}

$\kappa^s$ and $\kappa^h$  are  invariant under finite coverings. 
Thus $\kappa^s$ can be extended to virtually symplectic/complex manifolds, manifolds which 
are finitely covered by symplectic/complex manifolds:
If $X$ is finitely covered by a symplectic/complex  manifold $M$, then $v\kappa(X): =\kappa^{s/h}(M)$.

\medskip
\noindent {\bf $1-$dimensional fibers}

When $M=S^1\times Y$ it was conjectured by Taubes and confirmed by Friedl and Vidussi \cite{FV} that $M$ is symplectic if and only if
$Y$ is fibred.  In this case the additivity is verified in \cite{Z}. 

For general circle bundles  it is essentially understood in \cite{FV-bundle} which ones admit symplectic structures, and the virtual Betti number calculations in  \cite
{Baykur} provide ample evidence for the additivity of $\kappa^s$ for   $S^1-$bundles, ie.
$\kappa^s(M) =  \kappa^t(Y)$.

\medskip

\noindent {\bf $2-$dimensional fibers}

Thurston observed that surface bundles admit symplectic structures  if the fibers are homologically essential, and the converse is also true
   (\cite{W}). 
A consequence of the equality  $\kappa^l(g, h, n)=\kappa^s$  in \cite{DZ} is that $\kappa^s$ is additive for surface bundles over surfaces, ie.
$\kappa^s=\kappa^t(fiber)+\kappa^t(base)$.

\medskip

\noindent {\bf $3-$dimensional fibers}

It is not completely understood yet which mapping tori  admit symplectic structures. 
However, a  `virtual' progress has been  made.  
To put it in context recall that, by Agol's solution of virtual fibration conjecture, all irreducible 3-manifolds are virtually fibred except some graph manifolds.

Suppose $X$  fibers over the circle with fiber $Y$ 
and $Y$ is finitely covered by a mapping torus with fiber $F$. 
Applying  Luttinger surgery  to  $F\times T^2$, it is shown in \cite{LN} (cf. also \cite{BF}) that
\begin{enumerate}
\item If $g(F)=0$, then $X$ is virtually symplectic and $v\kappa(X)=-\infty$.

\item If $g(F)=1$, then $X$ is virtually symplectic if and only if $vb_1(X)\geq 2$.
Moreover, if $vb_1(X)\geq 2$, then $v\kappa(X)=0$. 

\item If $g(F)>1$, then $X$ is virtually symplectic with $v\kappa=1$.
\end{enumerate}
In particular, if $v\kappa^s(X)$ is defined then it satisfies the additivity
$v\kappa(X)=\kappa^t(F)+\kappa^t(Y)$. 
\medskip




\subsection{Behavior under Surgeries}

\subsubsection{Luttinger surgery} 

It is shown in \cite{HL} that   $\kappa^s$ is unchanged under Luttinger surgery along Lagrangian  tori (\cite{Luttinger}, \cite{ADK}). Combined with the diffeomorphism invariance of $\kappa^s$, this fact can be used
to distinguish non-diffeomorphic manifolds. 

In \cite{AP-08}, \cite{AP-10}, \cite{BK}, \cite{FdPR},  several symplectic manifolds homeomorphic
to  small rational manifolds  are constructed via Luttinger surgery. With $\kappa^s=2$ manifolds such as 
$\Sigma_g\times \Sigma_h$ with $g, h\geq 2$ as the building blocks, the invariance of $\kappa^s$ gives a quick alternative proof that these manifolds are small exotic manifolds. 
The invariance of $\kappa^s$ under Luttinger surgery  is also used effectively in \cite{AO} to construct
non-holomorphic Lefschetz fibrations with arbitrary $\pi_1$, starting from K\"ahler surfaces with
$\kappa^s=\kappa^s=1$.

\subsubsection{Genus 0 sum and rational blow-down}

 For  a genus zero sum, the computation $\kappa^s$  is largely carried out  in \cite{D} by Dorfmeister  and  finished in \cite{DLW}. 
The general behavior of $\kappa^s$  under a genus zero sum is non-decreasing. 
The rational blow-down of $-4$ spheres in $\C P^2\# 10\overline{\C P^2}$  is the most interesting case. The resulting manifolds include some Dolgachev surfaces, which have $\kappa^s=1$,
and the Enriques surface which has $\kappa^s=0$. 

General rational blow-down 
operations of Fintushel-Stern \cite{FS-rbd} and Symington \cite{Sy} have been used effectively for  the symplectic geography problem, and especially striking in  \cite{jP} and subsequent works by Fintushel, J. Park, 
Stern, Stipsicz and Szabo for  the exotic
geography problem. We postulate  that $\kappa^s$ also non-decreases  under an arbitrary rational blow-down. 
It will be nice to have a simple way to determine the change of $\kappa^s$
(as well as for the related  star surgery  in \cite{KS}).

\subsubsection{Positive genus sum}

Usher  systematically investigated $\kappa^s$ for positive genus  sum in \cite{U2}. 
His  calculation is interpreted in \cite{LY} in terms of properties of the adjoint class of
the gluing surfaces. 
It is further rephrased  in terms of 
relative Kodaira dimension of the summands in \cite{LZ} (cf. Section 5.1): If $(M, \omega)$ is 
the positive genus fiber sum of $(M_i, \omega_i)$ along $F_i\subset M_i$, then 
\begin{equation}\label{sum} \kappa^s(M, \omega)=\max\{\kappa^s(M_1, \omega_1, F_1), \kappa^s(M_2, \omega_2, F_2)\}.
\end{equation} 
Such a formula is especially effective in distinguishing exotic smooth structures on symplectic manifolds homeomorphic to small rational manifolds (eg. \cite{A-08}). 
Notice that this formula applies to Fintushel-Stern's powerful knot surgery (\cite{FS-link}), which  is a genus 1 fiber sum.


\section{Main problems and progress in each class}
In this section we will discuss properties of symplectic 4-manifolds in each $\kappa^s$  class. 
The slogan is: the smaller $\kappa^s$  the more we understand. 

\subsection{Surfaces and symmetry of  $\kappa=-\infty$ manifolds}

The symplectic 4-manifolds with Kodaira dimension $\kappa=-\infty$ have been
classified.  As smooth 4-manifolds, they are rational or ruled (ie. diffeomorphic to projective surfaces which are rational or ruled). And the moduli space
of symplectic structures is identified with the quotient of the symplectic cone by the geometric
automorphism  group. Moreover, both the symplectic cone and the geometric
automorphism group have been explicitly determined.

A central ingredient of  all the progress  is the existence and abundance of non-negative self-intersection symplectic surfaces.
  The major source of such surfaces is TaubesÕ symplectic Seiberg-Witten theory.
  Recall that the GT  (which means Gromov-Taubes)  invariant of a class $e$ in $H_2(M;\Z)$ is a Gromov type invariant defined by Taubes (cf. [61]) counting embedded (but not necessarily connected) symplectic surfaces representing the Poincar\'e dual to $e$, and $e$ is called a GT basic class if its GT invariant is nonzero.

There are still many open problems. Among them are: Symplectic versus K\"ahler, 
classifying symplectic surfaces and Lagrangian surfaces, determining the symplectomorphism group.  

\subsubsection{Smooth and symplectic classification}

The smooth classification  in \cite{Liu}   is achieved by finding 
a symplectic sphere with non-negative self-intersection. Such manifolds are uniruled, and a uniruled manifold is  a rational manifold or an irrational ruled manifold according to McDuff \cite{Mcjams}.

There are other characterizations in terms of smooth surfaces: 
\begin{itemize}

\item The existence of a smoothly embedded sphere with non-negative self-intersection (\cite{L-sphere}).

\item The existence of a smoothly embedded surface with positive genus $g$ and  self-intersection $2g-1$  (\cite{LMY}).

\end{itemize}
Notice that all $\kappa^s=-\infty$ manifolds have $b^+=1$. 
Due to the fundamental  fact that any $b^+=1$ manifold has infinitely many GT classes, 
 the inflation process of LaLonde-McDuff  (\cite{LMc-inflation}) can be applied  effectively, which is essential 
 for the symplectic classification (up to diffeomorphisms) of $\kappa^s=-\infty$ manifolds.


The moduli space  of symplectic structures, which is the space of diffeomorphic symplectic forms,  is completely understood (See \cite{Salamon} for a beautiful account).  

\begin{enumerate}

\item Symplectic structures are unique up to diffeomorphisms  in each cohomology class (\cite{LMc}, \cite{LLiu2}).

\item There is a unique symplectic canonical class up to orientation-preserving diffeomorphisms.

\item Therefore the moduli space $\mathcal M_X$ is identified with the quotient of the symplectic cone
by the geometric automorphism group. Here,  the symplectic cone is the open set of cohomology classes represented by symplectic forms,  and the  geometric automorphism group records the  action of the orientation-preserving diffeomorphism group on homology.

\item The symplectic cone   is completely determined in terms of the set $\mathcal E$ (\cite{LL}). 
And the subcone with  a fixed canonical class is a convex set, in particular, path connected.

\item The  geometric automorphism group  is generated by the reflections on $\mathcal E$,  $\mathcal L$ and $\mathcal H$, where $\mathcal L$ and $\mathcal H$ are the sets of 
the classes  represented by smoothly embedded spheres, and having square $-2$ and $1$ respectively
(\cite{FM}, \cite{LiL-generator}).

\item The sets $\mathcal E$, $\mathcal L$ and $\mathcal H$ are explicitly determined (\cite{LiL-genus}). 

\end{enumerate}

All $\kappa^s=-\infty$ manifolds admit K\"ahler structures. However, the following symplectic versus K\"ahler question is still open.   

\begin{question}  
For which $\kappa^s=-\infty$ manifold does there exist 
non-K\"ahler symplectic form?
\end{question}

 Surprisingly, it is shown in \cite{CP} that there are non-K\"ahler symplectic forms  on one point blow up of elliptic ruled surfaces. 
Due to properties  (1) and  (2), to show that any symplectic form is K\"ahler on a given
$\kappa^s=-\infty$ manifold, it suffices to show that the generic
K\"ahler cone coincides with the sub-symplectic cone with the fixed K\"ahler  canonical class. 
Such a equality of 
 cones  is known to be true for $S^2-$bundles (\cite{Mc-ruled}) and up to 9 blow-ups of the projective plane   (cf.  \cite{Bi}   for the connection to  the longstanding  Nagata conjecture and 
\cite{L-moduli} for results on other K\"ahler surfaces).

\subsubsection{Symplectic and Lagrangian surfaces}

The following is a conjecture  in \cite{DLW}  regarding  the existence of symplectic surfaces.

\begin{conj}\label{spec:isphere}Let $(M,\omega)$ be a symplectic $4$-manifold with $\kappa^s=-\infty$.
  Let $A\in H_2(M,\mathbb Z)$ be a homology class.
  Then $A$ is represented by a connected $\omega$-symplectic surface if and only if \begin{enumerate}
\item $[\omega]\cdot A>0$,
\item $g_\omega(A)\geq 0$, where $g_{\omega}(A)$ is defined by 
$2g_{\omega}(A)-2=K_{\omega}\cdot A+A\cdot A$,
\item $A$ is represented by a smooth connected surface of genus $g_{\omega}(A)$.
\end{enumerate}

\end{conj}

A rather general answer comes from Taubes' symplectic SW theory and the SW wall crossing formula
(\cite{KM}, \cite{LLiu1}): any  sufficiently large  $\omega-$positive class  is represented by
an $\omega-$symplectic surface.
In light of  properties (1), (2), and (4) of the moduli space of symplectic structures,   the following birational stability result in \cite{DLW} reduces the verification of the conjecture for an arbitrary symplectic form
to some fixed symplectic form $\omega_0$.

\begin{theorem}\label{stability} Let $(M,\omega)$ be a symplectic manifold with $b^+=1$ and
 $V$  an $\omega-$symplectic surface.
Then for any symplectic form $\tilde \omega$ homotopic to $\omega$ and positive on the class $[V]$, there exists an $\tilde \omega-$symplectic surface $\tilde V$, which  is smoothly isotopic to $V$.
\end{theorem}

This stability is applied  in \cite{DLW} to classify symplectic  spheres of negative self-intersection.

For orientable Lagrangian surfaces, the existence problem is completely solved. Since $\kappa^s=-\infty$ 
manifolds have $b^+=1$, 
there are no essential orientable Lagrangian surfaces  of positive genus.
For Lagrangian spheres there is the following simple criterion  in~\cite{LW} (compare with Conjecture \ref{spec:isphere}):
$A\in H_2(M;\mathbb Z)$  contains a Lagrangian sphere if and only if 
\begin{itemize}

\item $\omega(A)=0$,

\item  $A\cdot A=-2$,

\item  $A$ is represented by a smooth sphere.

\end{itemize}
In summary, a nonzero  class in $H_2(M;\Z)$  is represented by a Lagrangian surface if and only if it is in the set $\mathcal L$ and has zero pairing
with $\omega$. 
Furthermore, the solution has been  extended  to Lagrangian ADE configurations  in \cite{DLW}.
For uniqueness, it is shown in \cite{LW} and \cite{BLW} that homologous Lagrangian spheres are unique up to smooth isotopy. 
The much stronger uniqueness up to Hamiltonian isotopy was first established by Hind for the monotone $S^2\times S^2$ (\cite{Hind-Lag}) and then by Evans for monotone manifolds with Euler number at most $7$(\cite{Evans-Lag}). For further extensions see \cite{LW}. On the other hand, Seidel (\cite{Seidel}) discovered that Hamiltonian uniqueness fails for the monotone manifold with Euler number $8$.

There has been progress towards classifying Lagrangian $\R P^2$ in small rational manifolds (cf. \cite{BLW} for the status). With the recent classification of symplectic spheres with self-intersection $-4$ in \cite{DLW}, it seems possible  to classify Lagrangian $\R P^2$ since they correspond to symplectic $-4$ spheres
  via rational blow down. 


\subsubsection{Symplectic mapping class group and symplectic Cremona map}
A consequence of the simple criterion for the existence of Lagrangian spheres
is  the calculation of the homological action of the symplectomorphism group Symp$(M, \omega)$ in \cite{LW}: the action is generated by Dehn twists along Lagrangian spheres. 

The Torelli part is  much harder to determine. It is known to be   connected for rational manifolds with Euler number up to $7$, due to many people's work  (cf.  \cite{LLW} for  references). Hence the symplectic mapping class group 
is a finite reflection group  for these manifolds. 
On the other hand, Seidel showed that  the symplectic mapping class group for  the monotone 
rational manifold with Euler number $8$ is infinite (\cite{Seidel} and \cite{Evans-group}).  
The natural question is:  Is the symplectic mapping class group infinite for any symplectic rational manifold with Euler number at least $8$?

Finite symplectic symmetry is being investigated in  \cite{CLW}. This is the symplectic analogue of the classical Cremona problem in algebraic geometry. 
A classification of $\Z_n-$Hirzebruch surfaces 
up to orientation-preserving equivariant diffeomorphisms has been achieved in \cite{Chen-bundle}.

\subsection{Towards  the classification of $\kappa=0$ manifolds}


The basic problem is the  speculation on smooth classification.

\begin{conj}\label{conj-scy} (\cite{Do}, \cite{L06}) If $\kappa^s(M,\omega)=0$ and $(M, \omega)$ minimal,
then $M$ is diffeomorphic to one of the following: 

\begin{itemize}

\item K3, 

\item Enriques surface,   

\item a $T^2-$bundle over $T^2$

\end{itemize}

\end{conj}

Minimal $\kappa=0$ manifolds are exactly the ones with torsion $K_{\omega}$. 
If $(M, \omega)$ has $K_{\omega}=0$, $(M, \omega)$ is called  a symplectic Calabi-Yau surface.
In the case of SCY, the speculation is that  $M$ is diffeomorphic to K3 or a $T^2-$bundle over $T^2$.
Notice all the manifolds in the list  above allow some kind of torus fibrations.

\subsubsection{Homological classification}

Conjecture \ref{conj-scy} has been verified at the level of homology. 
The  crucial step is to derive the following bound on $b^+$.

\begin{theorem}\label{homology} (\cite{L}, \cite{Ba})
Suppose $\kappa(M,\omega)=0$. Then 
$b^+(M)\leq 3$. 

\end{theorem}

This  bound  on $b^+$ was established in \cite{MS} assuming  $b_1=0$, and in \cite{L06} assuming $b_1\leq 4$. 
The main idea is to apply
Furuta's $Pin(2)-$equivariant    finite dimensional  approximation (\cite{Furuta})  to the SW equation of the symplectic Spin$^c$ structure
to show that  $SW(K_{\omega})$ is even if $b^+>3$.
Then  invoke Taubes' fundamental calculation $SW(K_{\omega})=\pm 1$ (\cite{T0}).

Theorem \ref{homology} has the following consequences: If $M$ is minimal with $\kappa^s(M)=0$, then 

\begin{enumerate}

\item $b^-(M)\leq 19$. 

\item $b_1(M)\leq 4$.

\item The signature is equal to $0, -8, -16$. 

\item Euler number is equal to $0, 12, 24$.

\item $M$ either  has the same $\Z-$cohomology  group and intersection form as the K3 or the Enriques surface, or the  same $\Q-$homology group and intersection form  as  a $T^2-$bundle over $T^2$.

\item If $b_1(M)=4$, then $H^*(M;\Q)$ is generated by $H^1(M;\Q)$ and hence isomorphic to $H^*(T^4;\Q)$ as a ring. 
\end{enumerate}

Notice  that when $M$ is minimal, we have $$0=K_{\omega}\cdot K_{\omega}=2\chi(M)+3\sigma(M)=4-4b_1(M)+5b^+(M)+b^-(M),$$ from which 
both  bounds (1) and (2) follow. (3) then follows from the divisibility of $\sigma$ by $8$, and (4) is a consequence of $2\chi(M)+3\sigma(M)=0$.
The claim (5) is based on the Euler number bound (4) and the observation that a finite covering of a minimal $\kappa^s=0$ manifold is still such a manifold.
Finally, (6) relies on (5) and the main result in \cite{RS}.

\begin{remark}
A famous consequence of Yau's solution to the Calabi  conjecture is that any  K\"ahler manifold  with torsion canonical class  
admits Ricci flat metrics, and hence its $b_1$ is bounded by the real dimension (\cite{Yau}). Notice that this is still valid for  symplectic
$4-$manifolds with torsion canonical class due to consequence (2) above. 
Thus it was tempting to speculate the $b_1$ bound continues to hold for any symplectic manifold with torsion canonical class. 
However, Fine and Panov showed in \cite{FP-13} that this is far from true in dimension $6$ and higher (see also \cite{FP-09}, \cite{FP-10}). Fine-Panov's manifolds, which are obtained by Crepant resolutions of orbifold  twister spaces coming from hyperbolic geometry,  certainly cannot carry Ricci
flat metrics. Nevertheless, these singular twistor spaces  carry Chern-Ricci flat Almost K\"ahler metrics (cf. \cite{Pook}).
\end{remark}

The following table lists possible homological  invariants of minimal $\kappa=0$ manifolds:
$$\begin{tabular}{|c|c|c|c|c|c|}
\hline
$b_1$ & $b_2$ & $b^+$ & $\chi$ & $\sigma$ & known manifolds\\
\hline
0 & 22 & 3 & 24 & -16 & K3\\
\hline
0 & 10 & 1 & 12 & -8 & Enriques surface\\
\hline
4 & 6 & 3 & 0 & 0 & 4-torus\\
\hline
3 & 4 & 2 & 0 & 0 & $T^2-$bundles over $T^2$\\
\hline
2 & 2 & 1 & 0 & 0 & $T^2-$bundles over $T^2$\\
\hline
\end{tabular}
$$

Accordingly, there are three homology  types: K3 type, Enriques type, and torus bundle type, distinguished by the Euler number. 

All the known  $\kappa=0$ 
manifolds allow $T^2-$fibrations.
One approach towards the smooth classification is to detect the existence of   tori.
Suppose a homology K3 has a winding family,
then via parametrized SW theory 
there is an embedded symplectic torus for some symplectic form
in the winding family.
Existence of symplectic torus can be proved in   homology $T^2-$bundles with $T^2$. 

\subsubsection{Virtual $b_1$, SCY group, and partial  homeomorphism classification}
Recall that  $vb_1(M)$ is the supremum of $b_1(\tilde M)$ among all finite covers $\tilde M$. 
If $(M, \omega)$ has $\kappa^s=0$, then any finite cover  also has $\kappa^s=0$.  This implies that  $vb_1(M)\leq 4$. The following question is raised in \cite{FV-K=0}: 
Is the $\mathbb F_p-$virtual Betti number of a $\kappa^s=0$ manifold bounded by 4?

Any  minimal $\kappa^s=0$ symplectic $4-$manifold has torsion symplectic canonical class so  it admits
a finite cover with trivial symplectic canonical class. Thus such a manifold could be called a virtual SCY surface. 
Following \cite{FV-SCYgroup}, a finitely presented group $G$ is called a (v)SCY group if $G=\pi_1(M_G)$ for some (virtual) SCY surface $M_G$.
If $b_1(G)=0$  and $G$ is residually finite,  then $G=1$ or $\Z_2$ and the corresponding vSCY surfaces are  unique up to homeomorphism. 
\begin{itemize}
\item If $G=1$, then $M_G$ has the same intersection form as the K3 surface and hence  is homeomorphic to the K3 surface by Freedman's fundamental classification of simply connected topological 4-manifolds (\cite{MS}).

\item If $G=\Z_2$, then $M_G$ has the same intersection form and $w_2-$type as the Enriques surface and hence  is homeomorphic  to the Enriques surface by the extension in \cite{HK} 
of Freedman's classification to the case $\pi_1=\Z_2$ (\cite{L06}).

\end{itemize}
If $b_1(G)>0$ then  it follows from the consequence (5) of  Theorem \ref{homology} that
$$2\leq vb_1(G)\leq 4, \quad \chi(M_G)=\sigma(M_G)=0.$$ 
In this case, Friedl and Vidussi  showed in \cite{FV-SCYgroup}:
If $G=\pi_1$ of  a (Infra)solvable manifold,  then the corresponding SCY surfaces are unique  up to homeomorphism.

Since all $T^2-$bundles over $T^2$ are solvable manifolds,  the beautiful conclusion is that, 
 any known vSCY group $G$ determines the homeomorphism type of the corresponding vSCY surfaces.

\subsubsection{Constructions}
It has been verified that various constructions only produce minimal $\kappa^s=0$ manifolds in Conjecture \ref{conj-scy}. 
There are also a couple of constructions potentially producing new $\kappa=0$ manifolds. 

\medskip

\noindent{\bf Fiber bundles}

Suppose $(M, \omega)\to B$ is a fibre bundle and $\kappa^s(M)=0$. 

If $B$ is a surface  then $M$ is a $T^2$ bundle over $T^2$ (\cite{DZ}, \cite{Baykur}, \cite{FV-SCYgroup}). 

 For mapping tori with  prime fiber, 
 the fiber  has to be  a $T^2-$bundle over $S^1$, and $M$ is a $T^2-$bundle over $T^2$(\cite{LN}).

For {circle bundles},
the base  $B$ must be  a $T^2-$bundle (\cite{FV-K=0}).
In particular,  $M$ is a mapping torus with fiber $T^3$, and hence a $T^2-$bundle over $T^2$
by the claim above on mapping tori.

\medskip

\noindent{\bf Lefschetz fibrations/pencils}

Smith observed that  if  a SCY surface $M$  admits a genus $g$  Lefschetz fibration with singular fibers, then 
$g=1$ by the adjunction formula and 
it follows from  the classification of genus 1 Lefschetz fibrations of Moishezon+Matsumoto 
that $M$ is the K3 surface or a $T^2-$bundle over $T^2$ (\cite{S}). 
More generally, if   multiple fibers are allowed, then  $M$ is the Enriques surface.

On the other hand, according to Donaldson (\cite{Do-fibration}), any symplectic manifold admits Lefschetz pencils.
Baykur and collaborators are  able to 
  read off the Kodaira dimension from 
monodromy factorizations of Lefschetz pencils with multisections in
framed mapping class groups (\cite{BMV}). This could potentially lead to discovering new $\kappa^s=0$
manifolds distinguished by the fundamental group.

\medskip
\noindent{\bf Fiber sums}

If $(M, \omega)$ with $\kappa^s=0$  is a non-trivial  genus 0 fibred sum, then
$M$ is diffeomorphic to the blow up of the  Enriques surface, as  the rational blow down of $\C P^2\# l\overline{\C P}^2$ for some $l\geq 10$ (\cite{D}).

 If $(M, \omega)$ with $\kappa^s=0$  is a non-trivial positive genus fibred sum, then
the summands have $\kappa^s=-\infty$ (rational or ruled), and the sum  is along tori representing
$-K_{\omega}$ (\cite{U2}).
Moreover, if the sum is relatively minimal then $(M, \omega)$ is a vSCY surface by \cite{U1}.
In this case,  $M$ is diffeomorphic to the K3 surface, the Enriques surface,  or a $T^2-$bundle over
$T^2$ with $b_1=2$, and the only decompositions are as follows:
$K3=E(1)\#_f E(1)$ along fibers, 
Enriques =$E(1)\#_f (S^2\times T^2)$ along a fiber and a bi-section, 
and the sum of two $S^2-$bundles over $T^2$
along bi-sections give rise to $T^2-$bundles over $T^2$.
We remark that Akhmedov constructed   in \cite{Akhmedov-cy} 
a family of  simply connected symplectic CY 3-folds via fiber sums along $T^4$. 
\medskip

\noindent{\bf Hypersurfaces in $\C P^3$}

Smooth complex  hypersurfaces in $\C P^3$ of the same degree are diffeomorphic to each other, and  the degree 4 ones are just the K3 surface.

\begin{question} Is a degree 4 symplectic hypersurface in $\C P^3$ diffeomorphic to K3?
\end{question}

Such a symplectic hypersurface  is a symplectic CY surface
with $c_2=24$.
By the homological classification, $b_1$ vanishes and 
it is a homology K3.

\begin{remark} By the Lefschetz hyperplane theorem, any complex hypersurface is simply connected. 
This is known to be true for symplectic hypersurfaces up to degree 3. In fact, via pseudo holomorphic curve theory, 
Hind \cite{Hind} showed that degree 1 and 2 symplectic hypersurfaces are diffeomorphic to the complex hypersurfaces of the same degree.
Further, for degree up to $3$, a symplectic hypersurface 
has Kodaira dimension $-\infty$, thus the same conclusion for degree 3 can be shown to follow from the classification of $\kappa=-\infty$ manifolds. 
On the other hand, McLean (\cite{McLean}) recently constructed Donaldson type  symplectic hypersurfaces (\cite{Do-hypersurface}) with large degree and positive $b_1$. 

\end{remark}

\medskip

\noindent{\bf Luttinger surgery}

Luttinger surgery is a vSCY surgery and even preserves the K3, Enriques, $T^2-$bundle types (\cite{HL}). 
It is easy to check that some (possibly all) $T^2-$bundles over $T^2$ can be obtained from $T^4$ via Luttinger surgeries. 
More interestingly, if there is  a Lagrangian torus in the K3 such that the resulting manifold
after Luttinger surgery is not simply connected, then it would be
  a new  symplectic CY surface.  
We remark that the coisotropic surgery, which is a higher dimensional analogue of Luttinger surgery, has been explored in  \cite{Akhmedov-cy},  \cite{BK-cy} to construct 6-dimensional symplectic CY manifolds.
 
\subsection{Euler number and decomposition of $\kappa^s=1$ manifolds}
\subsubsection{On the non-negativity of   Euler number}

Gompf's family of manifolds in \cite{Gompf} which  any finitely presented group as the fundamental group have $\kappa^s=1$, so classification in
this case is unattainable. A fundamental conjecture is 

\begin{conj} \label{Euler}   If $\kappa^s(M)=1$ then its  Euler number $\chi(M)$ is non-negative.
\end{conj}

We have mentioned  that  the Euler number is indeed non-negative when $\kappa^s=0$.   In fact, the only known symplectic 4-manifolds with
negative Euler number are $g\geq 2$  $S^2-$bundles  and  their blow ups up to $4g-5$ points. 
Assuming the conjecture, the geography of $\kappa=1$ manifolds was investigated in \cite{BL}.

When $b^+=1$ this conjecture  holds. This implication is essentially contained in  \cite{Liu}, which says  that $b_1\leq 2$ if $b^+=1$ and $\kappa^s\geq 0$. If a  manifold with $b_1=2$ and $b^+=1$ 
has negative Euler number, then it must have $b_2=0$. But such a manifold has $\chi=-1$ and $\sigma=1$ and hence $K_{\omega}\cdot K_{\omega}=2\chi+3\sigma=1$. Since it has $b^-=0$ and hence minimal, such a manifold has $\kappa^s=2$. 
 
One  approach to this conjecture for manifolds with $b^+=2$ is to note that  
(i)  Gompf's manifolds are  constructed  via genus 1 fiber sum along square zero tori and (ii) the non-negativity of $\chi$ is preserved under a genus 1 fiber sum. 
Thus it is natural to postulate whether   manifolds with  $\kappa=1$ and $b^+=2$ can be split 
along symplectic tori into manifolds with $b^+=1$.
To find the splitting tori, a useful observation is that a genus 1 fiber sum often leads to non-trivial Gromov-Taubes counting for the resulting torus class (see \cite{Mc-GT}). 
When $b^+\geq 2$, there are only finitely many GT classes, and when $\kappa^s=1$ and minimal,  any GT class is a square zero class represented  by
a disjoint union of symplectic tori. 
Thus we are led to the simple (but hard) question: can we always find a splitting torus from these finitely many GT classes?

\subsubsection{ The GT length of $K_{\omega}$ and  conjectured upper bound on $b^+$}
When $b^+\geq 2$, the symplectic canonical class $K_{\omega}$ is  a GT class. 
 It is also conjectured in \cite{survey:Kod-dim4} that for a symplectic manifold with $\kappa^s=1$, 
\begin{equation} \label{Noether0} b^+\leq 3 +2\, l(K_{\omega}),\end{equation}
where $l(K_{\omega})$ is the GT length of $K_{\omega}$. 

\begin{definition}  A (fine) GT decomposition of a nonzero class $e$ is an unordered set of pairwise orthogonal nonzero GT  classes 
$\{A_1,\dots,A_m\}$ such that $e = A_1 + \cdots + A_m$. $m$ is called the length of the decomposition. The GT length $l(e)$ of the class $e$ is the maximal length among all such decompositions, and it is defined to be zero if  $e=0$ or $e$ is not a GT  class. 
\end{definition}

The inequality
\eqref{Noether0} is a variation of the  Noether type inequalities proposed in  \cite{MMP}, \cite{FS-canonical}.
Notice that  it  holds for elliptic surfaces $E(n)$, where $b^+=2n-1$ and $l(K_{\omega})=n-2$. It also seems not hard to verify that  \eqref{Noether0}  is preserved under a genus 1 sum. 
In addition, Theorem \ref{homology} can be interpreted as asserting  that the inequality \eqref{Noether0} holds  when $\kappa=0$ since  $l(K_{\omega})=0$ in this case.

\subsection{Geography and exotic geography of 
 $\kappa^s=2$ manifolds}

\subsubsection{Geography}
 
 The symplectic geography problem was  originally posed by Gompf in \cite{Gompf}.  
  In the case $\kappa^s=2$, it refers to the problem of determining which ordered pairs of {\it positive} integers are realized as $(\chi(M), \,\,K_{\omega}\cdot K_{\omega})   $ for some minimal symplectic 4-manifold $M$.
  There has been considerable progress: all the lattice points with 
 $K_{\omega}\cdot K_{\omega}-2\chi(M)=3\sigma(M)\leq 0$  have been filled 
 (see \cite{ABBKP} and references therein).
The region  of positive $\sigma$ is not well understood  yet (see  \cite{AHP} for the current status). 
 The basic  conjecture is 
 
 \begin{conj}\label{conj-bmy}  $\kappa^s=2$ manifolds satisfy the Bogomolov-Miyaoka-Yau inequality $K_{\omega}\cdot K_{\omega}\leq 3\chi(M)$.

\end{conj}

The BMY inequality  is valid for K\"ahler surfaces of general type, which is a classical theorem of  Miyaoka (\cite{Mi}) and Yau (\cite{Yau}).
LeBrun (\cite{LeB-Ein1}, \cite{LeB-Ein2})  verified this for 4-manifolds with   Einstein metrics, and 
recently Hamenstaedt (\cite{Ha}) verified it  for surface bundles over surfaces.  Symplectic manifolds near the BMY line have been constructed in \cite{Stip}, \cite{Stip-BMY} and \cite{APU}. 
An open problem is to construct symplectic non-K\"ahler  manifolds on the BMY line. 
We remark that Conjecture \ref{Euler} for minimal $\kappa^s=1$ manifolds  fits with Conjecture \ref{conj-bmy}. 

On the other  hand, there are  
also  the  Noether type inequalities proposed in  \cite{MMP}, \cite{FS-canonical}, which are  conjectured lower bounds for $K_{\omega}\cdot K_{\omega}$ (cf. a sightly stronger version in \cite{survey:Kod-dim4}, and 
 see \cite{FL} for progress).
 
 As in the $\kappa^s=1$ case, one approach to the  geography conjectures  is to decompose along tori:
 If $M$ has  $b^+\geq 2$ and contains a square zero class of symplectic torus with non-trivial Gromov-Taubes invariant,
can $M$ be split into a genus 1 fiber sum?
This would be the analogue of the toroidal decomposition in dimension 3 along $\pi_1-$injective tori. 
When $b^+\geq 2$,  we call a manifold atoroidal if there are no genus 1 GT invariants. The geography of such manifolds was investigated in \cite{FjPS}.

\subsubsection{Exotic geography}
It is understood now that most  positive  pairs have more than one, or infinitely many representatives as well. The current interest is on small pairs. Remarkable progress on  
exotic small manifolds (necessarily having $\kappa^s=2$)  has been made by  Akhmedov-D. Park,  J. Park, Fintushel-Stern, Stipsicz-Oszv\'ath-Szabo, Baldridge-Kirk  etc. (cf. discussions  in section 3.2).

There are at most countably many distinct smooth structures on a closed topological manifold, and it  has been suggested to use
\begin{itemize}
 \item the minimal genus function (\cite{Lawson},  \cite{DHL-genus}),
 \item size of the geometric automorphism groups (\cite{DHL-group}),
 \item size of  finite symmetry (\cite{CK})
 \end{itemize}
  to order
smooth/symplectic structures on a topological 4-manifold. The moral is that a smooth structure on a topological 4-manifold
is considered the `standard one' if  it has the smallest  minimal genus function, largest geometric automorphism group, or largest finite symmetry  among all smooth structures.  In \cite{DHL-genus} it is verified
that $S^2\times T^2$, $S^2-$bundles over $S^2$, the Enriques surface,   and $\C P^2 \# n\overline{ \C P^2 }$ with $ n\leq9$ are standard in the sense of having the smallest minimal genus
function.
It will be nice to show that  $\C P^2 \# n  \overline{\C P^2 } $ is `standard' for any $n$.

\section{Extensions}

\subsection {Relative Kodaira  dim for a symplectic pair}

 Let $(M, \omega)$
be  a connected, closed  symplectic $4-$manifold and $F\subset(M,
\omega)$  a symplectic surface, not necessarily connected but  having  no sphere components. For such  a pair the notion of   relative  Kodaira dimension is  introduced in \cite{LZ}.

The definition involves the  adjoint class of $F$, which is $K_{\omega}+[F]$.
Assume that $F$ is maximal in the sense that 
 its  adjoint class satisfies 
$$(K_{\omega}+[F])\cdot E\geq 0 $$
for any $E\in \mathcal E_{\omega}$. 
Clearly, this is the same as  $[F]\cdot E >0$ for any class $E$, or equivalently, the complement of $F$ is minimal. 
Thus we call a pair  minimal if   $F$ is maximal. 
As in the absolute case, any pair can be blown down to a minimal pair. What is a bit surprising is
that the minimal models are unique  in the relative setting.

The adjoint class of a maximal $F$ without sphere components  satisfies the following positivity. 

\begin{lemma}\label{adjoint}  Suppose $F$ is maximal and has no sphere components. 

If $\kappa^s(M, \omega)\geq 0$, then 
$$(K_{\omega}+[F])  \cdot [\omega]> 0, \quad (K_{\omega}+[F])^2\geq 0.$$

If $\kappa^s(M, \omega)=-\infty$ and   $(K_{\omega}+[F])^2>0$, 
then   $(K_{\omega}+[F])\cdot [\omega]>0$. 
\end{lemma}

In particular,  it is impossible to have a maximal surface without sphere components such that 
$$(K_{\omega}+[{F}])\cdot \omega=0\quad \hbox{ and}  \quad (K_{\omega}+[{F}])^2>0.$$
 

 \begin{definition}  \label{symp Kod}
Let $F\subset(M, \omega)$ be a  symplectic surface without
sphere components.
\begin{itemize}
\item If $F$ is empty, then 
$\kappa^s(M, \omega, F)$ is defined to be $\kappa^s(M, \omega)$.

\item Suppose $F$ is non-empty and maximal. Then 
\[
\kappa^s(M, \omega, F)=\left\{\begin{array}{cc}
-\infty & \hbox{if $(K_{\omega}+[{F}])\cdot \omega<0$ or\,\, $(K_{\omega}+[{F}])^2<0$},\\
0& \hbox{ if $(K_{\omega}+[{F}])\cdot \omega=0$ and $(K_{\omega}+[{F}])^2=0$},\\
1& \hbox{ if $(K_{\omega}+[{F}])\cdot \omega>0$ and $(K_{\omega}+[{F}])^2=0$},\\
2& \hbox{ if $(K_{\omega}+[{F}])\cdot \omega>0$ and $(K_{\omega}+[{F}])^2>0$}.\\
\end{array}\right.
\]

\item For a general pair,  the Kodaira dimension is defined to be that of its unique minimal model. 
\end{itemize}
\end{definition}

$\kappa^s(M, \omega, F)$ is well defined in light of Lemma \ref{adjoint}. Here are basic properties
of $\kappa^s(M, \omega, F)$:
\begin{enumerate}

\item Suppose $F_1, F_2\subset (M, \omega)$  are maximal symplectic surfaces without sphere components. If $[F_1] = [F_2]$, then $\kappa^s(M,  \omega, F_1)=\kappa^s(M, \omega, F_2)$.

\item $$\kappa^s(M,  \omega, F)\geq \kappa^s(M, \omega).$$

\item The formula \eqref{sum}  for a positive genus fiber sum holds.

\item Suppose a nonempty surface $F\subset (M, \omega)$ is maximal with each component of positive genus. Then

\begin{itemize}
\item $\kappa^s(M,\omega,F) = -\infty$
 if and only if $M$ is a genus $h\geq 1$ $S^2-$bundle, and $F$ is a section.
 
\item $\kappa^s(M, \omega, F ) = 0$ if and only if
$\kappa^s(M)=-\infty$  and $[F] = -K_{\omega}$.

\end{itemize}

\end{enumerate}

Although Lemma \ref{adjoint} is not valid for spheres, 
 $\kappa^s(M, \omega, F)$ can be extended  to an arbitrary  embedded symplectic surface  $F\subset(M, \omega)$ as follows:
Let $F^+$ be the surface obtained from $F$ by removing the sphere components, and 
  define  $\kappa^s(M, \omega, F)$ as $\kappa^s(M, \omega, F^+)$. 
   All the results still hold in this more general setting with obvious modifications.
We notice that the above definition is similar in one aspect to that of the Thurston norm of $3-$manifolds: the $2-$spheres have to be discarded. One explanation is that a $2-$sphere has $\kappa^t=-\infty$, so it behaves like the empty set.

Due to Property (1), it is possible to further extend $\kappa^s(M, \omega, F)$ to the case of $F$ being a symplectic surface with pseudo-holomorphic singularities, or a weighted symplectic surface.
 In algebraic geometry there is  the notion of the log Kodaira dimension of a noncomplete variety introduced by Iitaka (see \cite{Ib}). 
 \begin{question}  For a projective pair,  is the relative Kodaira dimension  equal to Iitaka' s log Kodaira dimension?
 \end{question}
The log Kodaira  dimension of cyclic affine surfaces has been extensively studied, so we could use them as testing ground. In fact, results  in  (\cite{McLean-log}) 
strongly suggest the symplectic nature of log Kodaira dimension for affine surfaces.

\subsection{Symplectic manifolds with concave boundary}

More generally, we  ask for what open symplectic manifolds  we can define the Kodaira dimension. 
An appropriate category consists of symplectic 4-manifolds with concave boundary. Such manifolds are also called symplectic caps. 
In \cite{LMY}, uniruled caps
and Calabi-Yau caps are introduced. 

\begin{definition}
 Let $(P,\omega)$ be a compact symplectic 4-manifold  with  concave boundary  $(Y,\xi)$.

$(P, \omega)$ is called a uniruled cap if 
 $[c_{1}(P)]\cdot[(\omega,\alpha)]>0$ for 
  a choice of Liouville contact one form $\alpha$ (induced by a choice of Liouville vector field $V$ defined near $Y$ pointing inward along $Y$).
 
 $(P, \omega)$ is called a Calabi-Yau cap if $c_1(P)$ is a torsion class. 
\end{definition} 

These  caps are useful to establish the finiteness of topological complexity of strong symplectic fillings and Stein fillings. In fact, all known contact 3-manifolds 
with some sort of bounded topological complexity of fillings   admit such caps.

We remark that since $[(\omega,\alpha)]$ is a relative class, $[c_{1}(P)]\cdot[(\omega,\alpha)]$ is well-defined.
The Kodaira dimension of a general  cap will be investigated in \cite{LM}. The uniruled caps are those  with Kodaira dimension $-\infty$, and
the Calabi-Yau caps are the minimal ones with Kodaira dimension $0$. 

 Suppose a cap $(P, \omega)$ is embedded in a closed manifold $(M, \Omega)$.
It is shown in \cite{LMY} that if  $(P, \omega)$ is uniruled then  
$\kappa^s(M, \Omega)=-\infty$, and  if $(P, \omega)$ is Calabi-Yau then $\kappa^s(M, \Omega)\leq 0$.  
We speculate  that $\kappa^s(P, \omega) \geq \kappa^s(M, \Omega)$ holds in general, which 
will be  the analogue of Property (2) for the relative Kodaira dimension.


\begin{thebibliography}{99} 
\bibitem{Akhmedov-cy} A. Akhmedov,  \textit{Symplectic Calabi-Yau 6-manifolds}, Adv. Math. 262 (2014), 115-125.
\bibitem{A-08} A. Akhmedov, \textit{Small exotic 4-manifolds}, Algebr. Geom. Topol. 8 (2008), 1781-1794


\bibitem{ABBKP}A. Akhmedov, S. Baldridge, I. Baykur, P. Kirk, and D. Park, \textit{Simply connected minimal symplectic
4-manifolds with signature less than $-1$}, J. Eur. Math. Soc. (JEMS), 12(1):133-161, 2010.
\bibitem{AHP} A. Akhmedov, M. Hughes, D. Park, \textit{Geography of simply connected nonspin symplectic 4-manifolds
with positive signature},  Pacific J. Math., 261(2):257-282, 2013.
\bibitem{AO} A. Akhmedov, B. Ozbagci, \textit{Exotic Stein Þllings with arbitrary fundamental group}, preprint, arXiv:1212.1743.

\bibitem{AP-08} A. Akhmedov, D. Park, \textit{Exotic smooth structures on small 4-manifolds}, Invent. Math. 173 (2008), no. 1, 209-223.
\bibitem{AP-10}A. Akhmedov, D. Park,  \textit{Exotic smooth structures on small 4-manifolds with odd signatures},  Invent. Math.,
181(3):577-603, 2010.
\bibitem{APU} A. Akhmedov, D. Park, G.Urz\'ua, \textit{Spin symplectic 4-manifolds near the Bogomolov-Miyaoka-Yau line}, JGGT, 4, (2010), 55-66. 

\bibitem{ADK} D. Auroux, S. Donaldson, L. Katzarkov, \textit{Luttinger surgery along Lagrangian tori and non-isotopy for singular symplectic plane curves},
Math. Ann. 326 (2003), no. 1, 185-203.

\bibitem{Ba} S. Bauer, \textit{Almost complex 4-manifolds with vanishing first Chern class}, J. Differential Geom., 79(1):25-32, 2008.
\bibitem{Baykur} I. Baykur, \textit{Virtual Betti numbers and the symplectic Kodaira dimension of fibred 4-manifolds},  arXiv:1210.6584.

\bibitem{Baykur-vb} I. Baykur,  \textit{Virtual betti numbers and the symplectic Kodaira dimension of fibered 4-manifolds}, 
Proceedings of the American Mathematical Society, 142 (2014), 4377-4384. 
\bibitem{BMV} I. Baykur, K. Hayano, \textit{Multisections of Lefschetz fibrations and topology of symplectic 4-manifolds}, arXiv:1309.2667.
\bibitem{BF} I. Baykur, S. Friedl, \textit{Virtually symplectic fibered 4-manifolds}, arXiv:1210.4983, to appear in   Indiana University Mathematics Journal.

\bibitem{BK} S. Baldridge, P. Kirk, \textit{Constructions of small symplectic 4-manifolds using Luttinger surgery}, J. Differential Geom. 82 (2009), no. 2, 317-361.
\bibitem{BK-cy} S. Baldridge, P. Kirk,   \textit{Coisotropic Luttinger surgery and some new examples of symplectic Calabi-Yau 6-manifolds}, 
 Indiana Univ. Math. J. 62 (2013), no. 5, 1457-1471.
\bibitem{BL}  S. Baldridge, T. J. Li,  \textit{Geography of symplectic 4-manifolds with Kodaira dimension one},  Algebr. Geom. Topol. 5 (2005), 355-368.
\bibitem{BLW} M. S.  Borman, T. J. Li,  W. Wu, \textit{Spherical Lagrangians via ball packings and symplectic cutting},
Selecta Math. (N.S.), 20(1):261-283, 2014.


\bibitem{Bi} P. Biran, \textit{From symplectic packings to algebraic geometry and back}, European Congress of Mathematics, Vol. II (Barcelona, 2000), 507-524, Progr. Math., 202, Birkhuser, Basel, 2001.
\bibitem{B} C. Birkar, \textit{Iitaka conjecture $C_{n,m}$ in dimension six}, arXiv:0806.4413
\bibitem{Chen-bundle} W. Chen, \textit{$G-$minimality and invariant negative spheres in $G-$Hirzebruch surfaces}, arXiv:1312.0848.
\bibitem{Chen-order} W. Chen, \textit{On the orders of periodic diffeomorphisms of 4-manifolds}, 
Duke Math. J. 156 (2011), no. 2, 273-310. 
\bibitem{CK} W. Chen, S. Kwasik, \textit{Symmetric symplectic homotopy K3 surfaces},
J. Topol. 4 (2011), no. 2, 406-430. 
\bibitem{CLW} W. Chen, T.J. Li, W. Wu, \textit{Symplectic Cremona map}, in preparation. 
\bibitem{CP} P. Cascini, Panov, \textit{Symplectic Generic Complex Structures on Four-Manifolds with $b^+=1$},  J. Symplectic Geom.
Volume 10, Number 4 (2012), 493-502.
\bibitem{DHL-group} B. Dai, C. Ho, T. J. Li, \textit{Geometric automorphism groups of symplectic 4-manifolds},  Topology Appl. 162 (2014), 1-11. 
\bibitem{DHL-genus} B. Dai, C. Ho, T. J. Li,   \textit{Minimal genus for 4-manifolds with $b^+=1$}, preprint. 

\bibitem{De} J.-P. Demailly, \textit{K\"hler manifolds and transcendental techniques in algebraic geometry}, International Congress of Mathematicians, Madrid 2006, Volume I, Plenary lectures and ceremonies, European Math. Soc. (2007) 153-186.



\bibitem{Do-hypersurface} S. Donaldson, \textit{Symplectic submanifolds and almost-complex geometry}, J. Differential Geom. 44
(1996), 666-705.

\bibitem{Do-fibration} S. Donaldson,  \textit{Lefschetz pencils on symplectic manifolds}, J. Differential Geom. 53 (1999), 205-236. 



\bibitem{Do} S. Donaldson,  \textit{Some problems in differential geometry and topology}, Nonlinearity 21 (2008) 157-164.

\bibitem{D-tori} J. Dorfmeister, \textit{Kodaira dimension of Lefschetz fibrations over tori}, to appear in Proc. AMS.

\bibitem{D}  J. Dorfmeister, \textit{Kodaira dimension of symplectic fiber sums along spheres}, arXiv:1008.4447,  Geom. Dedic. March 2014. 
\bibitem{DLW} J. Dorfmeister, T. J. Li,  W. Wu,    \textit{Stability and Existence of Surfaces in Symplectic 4-Manifolds with $b^+= 1$},
arXiv preprint arXiv:1407.1089.

\bibitem{DZ} J. G. Dorfmeister, W. Zhang, \textit{The Kodaira Dimension of Lefschetz Fibrations},  Asian J. Math. 13 (2009), no. 3, 341-357.
\bibitem{Evans-Lag} J. Evans, \textit{Lagrangian spheres in del Pezzo surfaces}, J. Topol. 3 (2010) 181-227
\bibitem{Evans-group} J. Evans, \textit{Symplectic mapping class groups of some Stein and rational surfaces}, J. Symplectic Geom. 9 (2011) 45-82. 
\bibitem{FL} P. Feehan, T. Leness, \textit{The SO(3) monopole cobordism and superconformal simple type}, arXiv:1408.5307.

\bibitem{FP-09} J. Fine, D. Panov,  \textit{Symplectic Calabi-Yau manifolds, minimal surfaces and the hyperbolic geometry of the conifold},  J. Differential Geom., 82 (1), 155-205, (2009).
\bibitem{FP-10} J. Fine, D. Panov, \textit{Hyperbolic geometry and non-K\"ahler manifolds with trivial canonical bundle},  Geometry and Topology 14, 1723-1763, (2010).
\bibitem{FP-13} J. Fine, D. Panov,      \textit{The diversity of symplectic Calabi-Yau 6-manifolds}, 
Journal of Topology, 2013.
\bibitem{FjPS} R. Fintushel, J. Park, R. Stern, \textit{Rational surfaces and symplectic 4-manifolds with one basic
class}, Algebr. Geom. Topol. 2 (2002), 391-402.
\bibitem{FS-rbd} R. Fintushel, R. Stern, \textit{Rational blowdowns of smooth 4-manifolds}, J. Diff.
Geom. 46, (1997), 181-235.


\bibitem{FS-link} R. Fintushel, R. Stern, \textit{Knots, links, and 4-manifolds}, Inv. Math. 134 (1998), 363-400.
\bibitem{FS-canonical} R. Fintushel, R. Stern, \textit{The canonical class of a symplectic 4-manifold}, Turkish J. Math. 25 (2001), 137-145. 
\bibitem{FdPR} R. Fintushel, B. Park, R. Stern, \textit{Reverse engineering small 4-manifolds}, Algebr. Geom. Topol. 7 (2007), 2103-2116.
\bibitem{FV} S. Friedl, S. Vidussi, \textit{Twisted Alexander polynomials detect fibered 3-manifolds},   Ann. of Math. 173 (2011).

\bibitem{FV-bundle}  S. Friedl, S. Vidussi, \textit{A vanishing theorem for twisted Alexander polynomials with applications to symplectic 4-manifolds}, J. Eur. Math. Soc. (JEMS) 15 (2013), no. 6, 2127-2041.

\bibitem{FV-K=0} S. Friedl, S. Vidussi, \textit{Symplectic 4-manifolds with K = 0 and the Lubotzky alternative}, Math. Res.
Lett. 18 (2011) 513-519.

\bibitem{FV-SCYgroup} S. Friedl, S. Vidussi, \textit{On the topology of symplectic Calabi-Yau 4-manifolds}, J. Topol. 6 (2013).


\bibitem{FM} R. Friedman, J. Morgan, \textit{Algebraic surfaces and
Seiberg-Witten invariants}, J. Alg. Geom., 6 (1997), 445-479.
\bibitem{FQ} R. Friedman, Z. Qin,
\textit{The smooth invariance of the Kodaira dimension of a complex surface},
Math. Res. Lett. 1 (1994), no. 3, 369-376. 
\bibitem{Furuta} M. Furuta,   \textit{Monopole equation and the 11/8 conjecture}, Math. Res. Letters 8 (2001),
279-291.

\bibitem{Gompf} R. Gompf, \textit{A new construction of symplectic manifolds}, Ann. Math., 142 (1995), 527-595.

\bibitem{Ha} U. Hamenst\"adt, \textit{Signatures of surface bundles and Milnor Wood inequalities},
arXiv:1206.0236.

\bibitem{HK} I. Hambleton, M. Kreck, \textit{Cancellation, elliptic surfaces and the topology of certain four-manifolds}, J. Reine Angew. Math. 444 (1993) 79Ð100. 

\bibitem{Hind-Lag} R. Hind, \textit{Lagrangian spheres in $S^2\times S^2$}, Geom. Funct. Anal, 14 (2004) 303-318.
\bibitem{Hind} R. Hind,    \textit{Symplectic hypersurfaces in $\C P^3$}, Proc. Amer. Math. Soc. 134 (2006), no. 4, 1205-1211. 


\bibitem{HL} C. Ho,  T. J. Li,  \textit{Luttinger surgery and Kodaira dimension},  Asian J. Math. 16 (2012), no. 2, 299-318.
\bibitem{I} S. Iitaka, \textit{On D-dimensions of algebraic varieties}, J. Math. Soc. Jap. 23(1971), 356-373
\bibitem{Ib} S. Iitaka, \textit{Algebraic geometry. An introduction to birational geometry of algebraic varieties}, Graduate Texts in Mathematics, 76. North-Holland Mathematical Library, 24. Springer-Verlag, New York-Berlin, 1982. x+357 pp.
\bibitem{IRS} M. Ishida, R. R\u{a}sdeaconu, I. Suvaina, \textit{On normalized Ricci flow and smooth structures on four-manifolds with $b^+=1$}, Archiv der Mathematik 92 (2009), no. 4, 355-365.
\bibitem{Ko} O. Kobayashi, \textit{Scalar curvature of a metric of unit
volume}, Math. Ann., 279 (1987), 253-265.
\bibitem{K} \textit{Kunihiko Kodaira: collected works. Vol. III.} Edited by Walter L. Baily, Jr.
Iwanami Shoten, Publishers, Tokyo; Princeton University Press, Princeton,
 N.J.,  1975. pp. i--x and 1142--1621.
 \bibitem{KM} P. Kronheimer,  T. Mrowka,  \textit{The genus of embedded surfaces in the projective plane}, Math. Res. Letters 1 (1994), 797-808. 
 \bibitem{KS} C. Karakurt, L. Starkston,    \textit{Surgery along star shaped plumbings and exotic smooth structures on 4-manifolds},  arXiv:1402.0801.  
\bibitem{LaL} G. LaNave,  T. J. Li, \textit{Symplectic Kodaira dimension via Donaldson's peak sections} in preparation. 
\bibitem{LMc} F. LaLonde, D. McDuff,  \textit{J-curves and the classification of rational and ruled symplectic 4-manifolds}, Contact and symplectic geometry Publ. Newton Inst. 8, (Cambridge,
1994), 3-42, Cambridge Univ. Press, Cambridge, 1996
\bibitem{LMc-inflation} F. LaLonde,  \textit{Isotopy of symplectic balls}, GromovÕs radius and the structure of ruled
symplectic 4-manifolds, Math. Ann. 300, no. 2, 273-296, (1994).
\bibitem{Lawson} T. Lawson,   \textit{The minimal Genus problem}, Expo. Math. 15(1997), 385-431. 

\bibitem{LeB-Ein1} C. LeBrun, \textit{Einstein metrics and Mostow rigidity}, MRL 2(1995), 1-8.

\bibitem{LeB-Ein2} C. LeBrun,  \textit{Einstein metrics, complex surfaces, and symplectic 4-manifolds}, MPCPS 147 (2009), 1-8.

\bibitem{LeB} C. LeBrun, \textit{On the Notion of general type},  Rend. Mat. Appl. (7)  17  (1997),  no. 3, 513--522
\bibitem{LeB1} C. LeBrun, \textit{Scalar curvature, covering spaces, and Seiberg-Witten theory},  New York J. Math.  9  (2003), 93--97
\bibitem{LeB2} C. LeBrun, \textit{Kodaira dimension and the Yamabe problem}, Comm. Anal. Geom.  7  (1999),  no. 1, 133--156.
\bibitem{L}T. J. Li, \textit{Quaternionic bundles and {B}etti numbers of symplectic
              4-manifolds with {K}odaira dimension zero}, Int. Math. Res. Not., 2006, Art. ID 37385, 28

\bibitem {L06} T. J. Li,  {\it Symplectic 4-manifolds with Kodaira dimension zero}, J. Differential Geometry 74 (2006) 321-352.

\bibitem{survey:Kod-dim4} T. J. Li, {\it  The Kodaira dimension of symplectic 4-manifolds}, Floer homology, gauge theory, and low-dimensional topology, 249Ð261, 
Clay Math. Proc., 5, Amer. Math. Soc., Providence, RI, 2006. 
\bibitem{L-moduli} T. J. Li, \textit{The space of symplectic structures on closed $4-$manifolds},    Third International Congress of Chinese Mathematicians. Part 1, 2, AMS/IP Stud. Adv. Math., 42, pt. 1, vol. 2, Amer. Math. Soc., Providence, RI, 2008, pp. 259Ð277.

\bibitem {L-sphere} T. J. Li, {\it Smoothly embedded spheres in symplectic
4-manifolds}, Proc. AMS. Vol. 127 No. 2 (1999) 609-613.
\bibitem{L-surface} T. J. Li, {\it Existence of symplectic surfaces},
Geometry and Topology of manifolds, 203-217, Fields Inst. Commun.,
47, AMS, 2005.
\bibitem{LiL-genus} B. H. Li, T. J. Li, {\it Symplectic genus, minimal genus and diffeomorphisms}, Asian J. Math. 6, (2002), 123-144.
\bibitem{LiL-generator} B. H.  Li, T. J. Li, \textit{On the diffeomorphism groups of rational and ruled 4-manifolds}, J. Math. Kyoto Univ. 46 (2006), no. 3, 583-593.
\bibitem{LLW} J. Li,  T. J. Li,  W. Wu, \textit{The symplectic mapping class group of $\C P^ 2 \#n {\overline {\C P^ 2}} $ with $n\leq 4$}, arXiv:1310.7329, to appear in Michigan Math. J. 
\bibitem{LL} T. J. Li, A. Liu, \textit{Uniqueness of symplectic canonical class, surface cone and symplectic cone of 4-manifolds with $B^+=1$},  J. Differential Geom.  58  (2001),  no. 2, 331--370.
\bibitem {LLiu1} T. J. Li, A. Liu, {\it Symplectic
structures on ruled surfaces and a generalized adjunction
inequality}, Math. Res. Letters 2 (1995), 453-471.
\bibitem{LLiu2} T. J. Li, A. Liu, \textit{General wall crossing formula},
Mathematical Research Letters 2, 797-810.
\bibitem{LMY} T. J. Li,  C. Y. Mak, \textit{Uniruled caps and Calabi-Yau caps},   arXiv: 1412.3208.
\bibitem{LM} T. J. Li, C. Y.  Mak, \textit{Kodaira dimension of symplectic 4-manifolds with contact boundary}, in preparation.
\bibitem{LN} T. J. Li,  Y. Ni, \textit{Virtual Betti numbers and virtual symplecticity of 4-dimensional mapping tori},  Math. Z. 277 (2014), no. 1-2, 195-208.
\bibitem{LR} T. J. Li, Y. Ruan, \textit{Symplectic Birational Geometry}, 307-326, { New Perspectives and Challenges in Symplectic Field Theory}, AMS. 2009.  
\bibitem{LW} T. J. Li,  W. Wu,  \textit{Lagrangian spheres, symplectic surfaces and the symplectic mapping class group}, Geom. Topol. 16 (2012), no. 2, 1121-1169.
\bibitem{LZ} T. J. Li, W. Zhang, \textit{Additivity and relative Kodaira dimensions},  Geometry and analysis. No. 2, 103Ð135, Adv. Lect. Math. (ALM), 18, Int. Press, Somerville, MA, 2011.
\bibitem{Liu} A. Liu, \textit{Some new applications of the general wall crossing formula}, Math. Res. Lett. 3 (1996) 569-585.
\bibitem{LY} T. J. Li, S.T. Yau, \textit{Embedded surfaces and Kodaira dimension}, preprint.

\bibitem{Luttinger}K. M.  Luttinger, \textit{Lagrangian tori in $\R^4$}, J. Diff. Geom. 42, 220-228 (1995).
\bibitem{MMP} M. Marino, G. Moore, G. Peradze, \textit{Four-manifold geography and superconformal symmetry},
Math. Res. Lett. 6 (1999), 429-437.
\bibitem{Mcjams} D. McDuff, \textit{The structure of rational and ruled symplectic
$4-$manifolds}, J. Amer. Math. Soc.  3 (1990), 679-712.
\bibitem{Mc-ruled}D. McDuff, \textit{Notes on ruled symplectic 4-manifolds}, Trans. Amer. Math. Soc. 345
(1994), no. 2, 623-639.
\bibitem{Mc-isotopy} D. McDuff, \textit{From symplectic deformation to isotopy}, (Irvine, CA, 1996), 85-99, First Int. Press Lect Ser. I, Internat. Press, Cambridge, MA, 1998.\bibitem{Mc-GT} D. McDuff, {\it Lectures on Gromov invariants for symplectic 4-manifolds}, Gauge theory and symplectic geometry, 1997, 175-210.
\bibitem{McS} D. McDuff, D. Salamon, \textit{A survey of symplectic $4$-manifolds with $b_+=1$},  Turkish J. Math.  20  (1996),  no. 1, 47--60.
\bibitem{McLean} M. McLean,  Private communication, 2014. 
\bibitem{McLean-log} M. McLean,    \textit{Symplectic invariance of uniruled affine varieties and log Kodaira dimension}, 
Duke Math. J. 163 (2014), no. 10, 1929-1964.
\bibitem{Mi} Y. Miyaoka, \textit{On the Chern numbers of surfaces of general type}, Inv. Math. 42  (1977), 225-237.
\bibitem{MS} J. Morgan, Z. Szabo,  \textit{Homotopy K3 surfaces and Mod 2 Seiberg-Witten invariants}, Math. Res. Lett. 4 (1997) 17-21.

\bibitem{OO} H. Ohta, K. Ono,  \textit{Notes on symplectic 4-manifolds with $b_2^+=1$, II}, IJM 7 (1996), 755-770.
\bibitem{jP} J. Park, \textit{Simply connected symplectic 4-manifolds with $b^+_2=1$ and $c_1^2=2$},
Invent. Math. 159 (2005), no. 3, 657-667. 
\bibitem{Pook} J. Pook, \textit{Homogeneous and locally homogeneous solutions to curvature flow},
arXiv: 1202.1427.
\bibitem{R} R. R\u{a}sdeaconu, \textit{The Kodaira dimension of diffeomorphic K\"ahler threefolds}, Proc. Amer. Math. Soc. 134 (2006), no.12, 3543--3553.
\bibitem{RS} D. Ruberman, S. Strule,\textit{Mod 2 Seiberg-Witten invariants of homology tori}, Math. Res. Lett. 7 (2000), 789-799.
\bibitem{Salamon} D. Salamon, \textit{Uniqueness of symplectic structures}, Acta Math. Viet. 38, (2013), 123-144. 

\bibitem{Seidel} P. Seidel, \textit{Lectures on four-dimensional Dehn twists}, from: 
"Symplectic 4-manifolds and algebraic surfaces", (F Catanese, G Tian, editors), Lecture Notes in Math. 1938, Springer, Berlin (2008) 231-267
\bibitem{S} I. Smith, \textit{Torus fibrations on symplectic four-manifolds}, Turkish J. Math. 25 (2001), no. 1, 69--95.
\bibitem{S2} I. Smith, \textit{Lefschetz pencils and divisors in moduli space}, Geom. Topol. 5 (2001), 579--608
\bibitem{Stip} A. Stipsicz, \textit{Simply connected symplectic 4-manifolds with positive signature}, In Proceedings of 6th G\"okova Geometry-Topology Conference, volume 23, pages 145-150, 1999.
\bibitem{Stip-BMY}A. Stipsicz, \textit{Simply connected 4-manifolds near the Bogomolov-Miyaoka-Yau line}, Math. Res. Lett., 5(6):723-730,
1998.
\bibitem{Su} I.  Suvaina, \textit{Yamabe invariant of a class of symplectic manifolds}, 
 arXiv:1410.1236.
 \bibitem{Sy} M. Symington,  \textit{Symplectic rational blowdowns}, J. Differ. Geom. 50, 505-518 (1998).
 
\bibitem{T0} C. H. Taubes,  \textit{The Seiberg-Witten invariants and symplectic forms}, Math. Research Letters, 1 (1994) 809-822.
\bibitem{T} C. H. Taubes, \textit{SW $\Rightarrow$ Gr: From the Seiberg-Witten equations to pseudoholomorphic curves}, J. Amer. Math. Soc. 9(1996), 845-918
\bibitem{T1} C. H. Taubes, \textit{Counting pseudoholomorphic
submanifolds in dimension $4$}, Jour. Diff. Geom. 44 (1996),
818-893.
\bibitem{Torres}     P.  Su‡rez-Serrato, R.  Torres, \textit{
A note on collapse, entropy, and vanishing of the Yamabe invariant of symplectic 4-manifolds},  arXiv:1408.1586, to appear in JGP. 

\bibitem{Ue} K. Ueno, \textit{Classification theory of algebraic varieties and compact complex spaces},
Notes written in collaboration with P. Cherenack. Lecture Notes in Mathematics, Vol. 439. Springer-Verlag, Berlin-New York, 1975. xix+278 pp.
\bibitem{U1} M. Usher, \textit{Minimality and symplectic sums},  Int. Math. Res. Not.  2006, Art. ID 49857, 17 pp.
\bibitem{U2} M. Usher, \textit{Kodaira dimension and symplectic sums}, To appear in Comment. Math. Helv.

\bibitem{W} R.  Walczak, \textit{Existence of symplectic structures on torus bundles over surfaces},
Ann. Global Anal. Geom. 28 (2005), no. 3, 211-231.
\bibitem{Yau} S. T. Yau, \textit{Calabi conjecture and some new results in algebraic geometry}, PNAU, 74, 1977, 1789-1799

\bibitem{Z}  W. Zhang, \textit{Geometric structures, Gromov norm and Kodaira dimensions},
arXiv:1404.4231.


\end{thebibliography}
\end{document}